\documentclass[12pt]{article}
\usepackage{amsfonts}
\usepackage{amsmath}
\usepackage{amssymb}
\usepackage{graphicx}

\usepackage{amsfonts}
\usepackage{amsmath}
\newtheorem{theorem}{Theorem}[section]
\newtheorem{lemma}[theorem]{Lemma}
\newtheorem{corollary}[theorem]{Corollary}

\newcommand{\IZ}{\mathbb Z}
\newcommand{\IC}{\mathbb C}
\newcommand{\IH}{\mathbb H}

\def\IH{{\Bbb H}}

\def\IS{{\Bbb S}}

\def\IZ{{\Bbb Z}} 
 
\def\IC{\Bbb C} 

\def\tr{{\rm tr}} 
\def\oC{\hat{\IC}} 
\def\tr{\mbox{\rm{tr\,}}}

\def\arccosh{\mbox{\rm{arccosh}}}

\def\ax{\mbox{\rm{ax}}}

\def\log{\mbox{\rm{log}}}

\def\em{\it}

\title{The continuous part of the axial distance spectrum for Kleinian groups.}
\author{G.J. Martin\thanks{This research was supported in part by a grant from the NZ Marsden Fund.\newline
\newline
{\bf  Mathematics subject classification 2010,}{Primary: 30F40, 57M10. } \newline 
{\bf  Keywords,} {Kleinian group, universal constraints, hyperbolic geometry}
}}
\date{}
\begin{document}

\maketitle

\begin{abstract}
Elements $f$ of finite order in the isometry group of hyperbolic three-space $\IH^3$ have a hyperbolic line as a fixed point set,  this line is the axis of $f$.  The possible  hyperbolic distances between axes of elements of order $p$ and $q$,  not both two, among {\em all}  discrete subgroups $\Gamma$ of  $Isom^+(\IH^3)$ has an initial discrete spectrum 
\[ 0 =\delta_0< \delta_1 < \delta_2 < \ldots <\delta_\infty,\]
each value taken with finite multiplicity, and above $\delta_\infty$ this spectrum of possible distances is continuous.  The value $\delta_\infty$ is the smallest number with the property that for each $\lambda<1$ there are only finitely many discrete groups generated by elements of order $p$ and $q$ whose axes are no more than $\lambda \delta_\infty(p,q)$ apart.  Geometrically $\delta_\infty$  places a bound on embedded tubular neighbourhoods of components of the singular set in the orbifold quotients $\IH^3/\Gamma$ and provides other geometric information about   this set.  The value $\delta_1(p,q)$ is known and tends to $\infty$ with $\min\{p,q\}$. Here we seek to determine - actually find asymptotically sharp upper-bounds for - $\delta_\infty(p,q)$.  We also show that  the gap $\delta_\infty(p,q)-\delta_1(p,q)$ is surprisingly small,  less than $1.4059\ldots$,  the sharp value for the Fuchsian case,  independent of $p$ and $q$.  This is despite both of these numbers tending to $\infty$ with either $p$ or $q$.
\end{abstract}
\newpage

\section{Introduction}

A {\em Kleinian group} $\Gamma$ is a discrete subgroup of the orientation preserving isometry group of hyperbolic $3$-space $Isom^+(\IH^3)$ which does not have an abelian subgroup of finite index.  We will use the identification $PSL(2,\IC)\cong Isom^+(\IH^3)$. Discreteness implies that $\Gamma$ acts properly discontinuously on $\IH^3$ and the orbit space $\IH^3/\Gamma$ is a hyperbolic $3$-orbifold. In fact noting that $\oC=\partial \IH^3/\Gamma$ and the isomorphism
\[ PSL(2,\IC)\ni  \left(\begin{array}{cc} a & b \\ c &d \end{array} \right) \leftrightarrow \frac{az+b}{cz+d} \in Con(\oC), \quad\quad ad-bc=1 \]
where $Con(\oC)$ is the group of conformal homeomorphisms of the Riemann sphere, the isomorphism $PSL(2,\IC)\cong Isom^+(\IH^3)$ is then effected by the Poincar\'e extension \cite[\S 3.3]{Beardon}.   

The elements   $f\in \Gamma$,  $f\neq 1$, fall into one of three classes determined by the value of
\begin{equation}
\beta(f)=\tr^2(f)-4 = 4\sinh^2 \frac{\tau+i\theta}{2}
\end{equation}
\begin{enumerate}
\item $f$ is elliptic: $\beta\in [-4,0)$ and $f$ is conjugate to a periodic rotation $z\mapsto \zeta z$, $\zeta$ a root of unity,  $\tau=0$, $\theta=\frac{p\pi}{q}$.
\item $f$ is parabolic: $\beta=0$ and $f$ is conjugate to the translation $z\mapsto z+1$.
\item $f$ is loxodromic:  $\beta \in \IC\setminus [-4,0]$, and $f$ is conjugate to the scaling $z\mapsto \lambda  z$, $|\lambda|\neq1$,  $\tau>0$.
\end{enumerate}
Thus if $f$ is elliptic or loxodromic,  it fixes two points and therefore the hyperbolic line in $\IH^3$ joining these points.  This line is call the {\em axis of} $f$,  denoted $\ax(f)$.
 \medskip 
 
 If $f$ is elliptic of order $p\geq 7$ in a discrete group $\Gamma$, then the axis of $f$ is simple; for each $g\in\Gamma$ either $g(\ax(f))=\ax(f)$ or $g(\ax(f))\cap\ax(f)=\emptyset$.  This follows from the classification of the elementary discrete groups, \cite[\S 5.1]{Beardon}.  If $f$ is a loxodromic element of short translation length,  then it also has a simple axis.  This is a well known consequence of J\o rgensen's inequality \cite{Jorgensen} and the formula (\ref{cosh2r}) below.   
 
 \medskip
 
 More is true for elliptic elements in a discrete group.  If $f$ is an elliptic element of order $p$ and $g$ is elliptic element of order $q$ in a discrete group $\Gamma$ and if $f$ and $g$ do not have the same axis,  then the axes of $f$ and $g$ are disjoint unless $\langle f,g\rangle$ is a finite spherical triangle group and $p,q\in \{2,3,4,5\}$, a finite dihedral group and $p$ or $q$ is equal to $2$, or a Euclidean group $p,q\in\{2,3,4,6\}$.  
 
 \medskip
 
 If the axes of elliptics $f$ and $g$ are sufficiently far apart the the group $\langle f,g\rangle$ is discrete and freely generated.  In \cite{GMM} we show the following.
 \begin{theorem}\label{free} Let $f,g$ be elliptics of order $p,q$,  not both $2$, respectively and suppose that the distance between the axes of $f$ and $g$, $\rho(f,g)$, satisfies
 \begin{equation}\label{6}
 \rho(f,g) \geq  \arccosh \left(\frac{1+\cos\frac{\pi}{p} \cos\frac{\pi}{p}}{\sin\frac{\pi}{p} \sin\frac{\pi}{p}}\right).
 \end{equation}
 Then $\langle f,g\rangle\subset Isom^+(\IH^3)$ is discrete and freely generated by $f$ and $g$.  This bound is sharp and achieved in the $(p,q,\infty)$ triangle group;  for $\frac{1}{r}<1-\frac{1}{p}-\frac{1}{q}$ the distance $\delta$ between the axes of the generating elliptics of orders $p$ and $q$ in the $(p,q,r)$-triangle group satisfy 
\begin{equation}\label{triangle} \delta =\arccosh \left(\frac{\cos\frac{\pi}{r}+\cos\frac{\pi}{p} \cos\frac{\pi}{p}}{\sin\frac{\pi}{p} \sin\frac{\pi}{p}}\right) \end{equation}
This triangle group is not freely generated and as $r\to\infty$ approaches the bound at (\ref{6}).
 \end{theorem}
 Conversely we record the following combined results from \cite{GMProc,GM1,GMM2}.
\begin{theorem} \label{delta1} Let $\delta_1(p,q)$ be the infimum of the non-zero distances between the axes of elliptics of order $p$ and $q$ in a discrete group. Then this infimum is attained unless $p=q=2$ and for $p\geq 7$, $q\geq 2$,
\[2\sin\frac{\pi}{p}\sin\frac{\pi}{q} \cosh \delta_1(p,q) 
= \left\{\begin{array}{ll} 1 & q\neq 3,p \\  \\ \cos\frac{\pi}{p} & q=3 \\ \\
\cos\frac{2\pi}{p} & q=p
\end{array}  
\right.\]
For $2\leq p \leq q \leq 6$ we have the table

\medskip
\begin{center}
{\bf Table 1. Values of  $2\sin\frac{\pi}{p}\sin\frac{\pi}{q} \cosh \delta_1(p,q) $}\\
\begin{tabular}{|c|c|c|c|c|c|}
\hline
& $2$ & $3 $& $4$ & $5$ & $6$   \\
 \hline
$2$ & & $1.765$ & $1.538$ & $1.300$ &$1.224$  \\
 \hline
$3$ & $1.765$ &$1.618$ & $\sqrt{2}$ & $1.401$ & $1$   \\   
\hline
$4$ & $1.538$ & $\sqrt{2}$ & $1.366$ & $1$  & $1$ \\
\hline
$5$ & $1.300$ &$1.401 $& $1$&$1$&$1$ \\
\hline
$6 $& $1.224$ & $1$ & $1$ & $1$& $1$  \\
\hline 
\end{tabular}
\end{center}
Each of these values is obtained by  either  a subgroup  of a  Coxeter group, that is a group generated by reflections, or a $\IZ_2$ extension of such a group.  In particular
\begin{enumerate}
\item For $p,q\leq 6$, the examples are related to tetrahedral reflection groups - groups generated  by reflections in the faces of a  hyperbolic tetrahedron. See {\rm \cite{GM1}.}
\item For  $p\geq 7$ and $q  =  2, 3, p$,  the  examples  are  the  triangle  groups  with
signature $(2, 3, p)$. These are the only examples  which are Fuchsian.
\item For $p \geq 7$ and $q = 4, 5$,  the examples are pentahedral reflection groups. See {\rm\cite{CM}.}
\item For $p,q\geq 7$ , $p\neq q$,   the examples are  the  hexahedral  reflection  groups. See {\rm \cite[\S 3]{GMM2} }
\end{enumerate} 
The  extremals  are  all unique  except  when  $p\geq 7$  and  $q  = 2$ when there are two examples.
\end{theorem}
 
 \medskip
 
 Our interest in this article is to supplement this information with further data of the possible distances between elliptic axes in a discrete group.  It turns out that there is initially a discrete set of possible values - the smallest possible value is given in Theorem \ref{delta1} - and a region beyond which any value is possible - Theorem \ref{free}.
 In particular,  here we seek the smallest value beyond which this spectrum of allowable axial distances is continuous. This is not the same as Theorem \ref{free},  there are {\em always} finite co-volume lattices in the continuous spectrum - indeed infinitely many. Our results are described in \S \ref{mainresult} and \S 5.
 
 \medskip
 
The base of the continuous spectrum puts a limit on possible {\em collaring theorems} which provide lower bounds on the distance between the axes of elliptic elements of $\Gamma$.  For instance the distance between the axis of $f$ and its nearest translate $g(\ax(f))=\ax(gfg^{-1})$. Such bounds give quantitative lower bounds on the radii of embedded tubular neighbourhoods of  the singular set,  a possibly trivalent graph formed by the projection of the elliptic fixed points,  in the three-orbifold $\IH^3/\Gamma$.  This  in turn provides lower bounds on their volume. If an elliptic axis is not simple,  then the group $\Gamma$ contains either a parabolic element or a spherical triangle group as a subgroup.  The elliptic axes  emanating from a finite point $x_0\in \IH^3$ fixed by a spherical triangle group are in general position and from this,  collaring bounds can be used to bound from below the distance from $x_0$ to its nearest translate,  and in turn provide lower bounds on the co-volume of $\Gamma$.

Remarkably the first several values of the spectrum of possible radii of maximal collars about elliptic axes of low order are realised only in arithmetic lattices.  See \cite{ADV}.    If the elliptic has high order,  then it has a large smallest maximal collar which provides enough volume to exceed that of known small co-volume arithmetic lattices.  These ideas led to the identification of the smallest co-volume lattices of Hyperbolic $3$-space, \cite{GMannals,MMannals}, a problem Siegel raised in 1942, and the classification programme for arithmetic Kleinian groups generated by torsion, \cite{MM,GMM2}.

\section{The axial distance spectrum for elliptics.}

Given two non-parabolic elements $f,g\in \Gamma$,   we define
\begin{equation}
\rho(f,g)=\rho_{\IH^3}(\ax(f),\ax(g))
\end{equation}
the hyperbolic distance between the two axes. We set $\delta_0(p,q)=0$, and inductively define for $f$ elliptic of order $p$ and $g$ elliptic of order $q$
\begin{equation}
\delta_1(p,q) = \inf\{\rho(f,g): \langle f,g \rangle \mbox{ is discrete and $\rho(f,g)>0$,} \} 
\end{equation}
\[
\delta_n(p,q) = \inf\{\rho(f,g): \langle f,g \rangle \mbox{ is discrete and $\rho(f,g) > \delta_{n-1}(p,q)$}, \} 
\]
and 
\[
\delta_\infty(p,q) = \inf\{\rho(f,g): \langle f,g \rangle \mbox{ is discrete and freely generated by$f$ and $g$}. \} 
\]

Our first result gives  some general facts about this spectrum which is non-decreasing by definition. 

\begin{theorem}\label{ads}  For all $n\geq 0$,  each value $\delta_n(p,q)$ is achieved in a Kleinian group and with finite multiplicity. The set $\{\delta_n(p,q)\}\cup\{\delta_\infty(p,q)\}$ is closed and bounded, 
\begin{equation}\label{ub}
\sup \{ \delta_n(p,q) \}  \leq \delta_\infty(p,q) < \frac{1+\cos \frac{\pi}{p}\cos \frac{\pi}{q}}{\sin \frac{\pi}{p}\sin\frac{\pi}{q}}.
\end{equation}
For each $p$ and $q$
\begin{equation}
\delta_n(2,p)\leq \delta_{n}(p,q), \quad\quad \delta_n(2,q)\leq \delta_n(p,q).
\end{equation}

For each $\delta<\delta_\infty(p,q)$,  there are only finitely many Kleinian groups generated by elliptics of orders $p$ and $q$ whose axes are no more than $\delta$ apart.

For each $\delta>\delta_\infty(p,q)$,  there are infinitely many Kleinian groups generated by elliptics of orders $p$ and $q$ whose axes are  $\delta$ apart.
\end{theorem}

\noindent{\bf Remark.}  It seems a challenging question to decide if the set $\{\delta_n(p,q)\}_{n\geq1}$  is finite or infinite. Certainly there are infinitely many possible values for the distance between axes.  Essentially the question  becomes whether or not  there are infinitely many values smaller than the smallest value for a group which  free on its generators.   The ellipse
\[ \{z = 4 \sin^2\frac{\pi}{p}\sin^2\frac{\pi}{q} \sinh^2(\delta_\infty(p,q)+i\theta):\theta\in [0,2\pi] \} \]
bounds a region in $\IC$,  call it ${\cal E}_{p,q}$,  which lies in  the complement of a moduli space of discrete groups freely generated by two elliptics of orders $p$ and $q$. This latter set is connected - it is the Schottky representation of the moduli space of the two-sphere with four cone points,  two of order $p$ and two of order $q$,  $\IS^2_{(p,p,q,q)}$.  This moduli space  is topologically an annulus as it is determined by the gluingings of a rigid pair of disks with cone points $p,p$ in one and $q,q$ in the other. The two parameters being the length of the boundary of the disk and the gluing angle.  This is not entirely straightforward,  but the details would take us too far astray. Within ${\cal E}_{p,q}$ are commutator values $\tr[f,g]-2$ determining discrete groups which are not free on their generators (where we view $f,g\in PSL(2,\IC)$).
This is illustrated below in Figure 1 for $p=2$, $q=3$  with data from nearly 1500 orbifold Dehn surgeries on two-bridge knot and link complements from \cite{Qing}.
\begin{center}
\scalebox{0.65}{\includegraphics*[viewport=20 460 660 660]{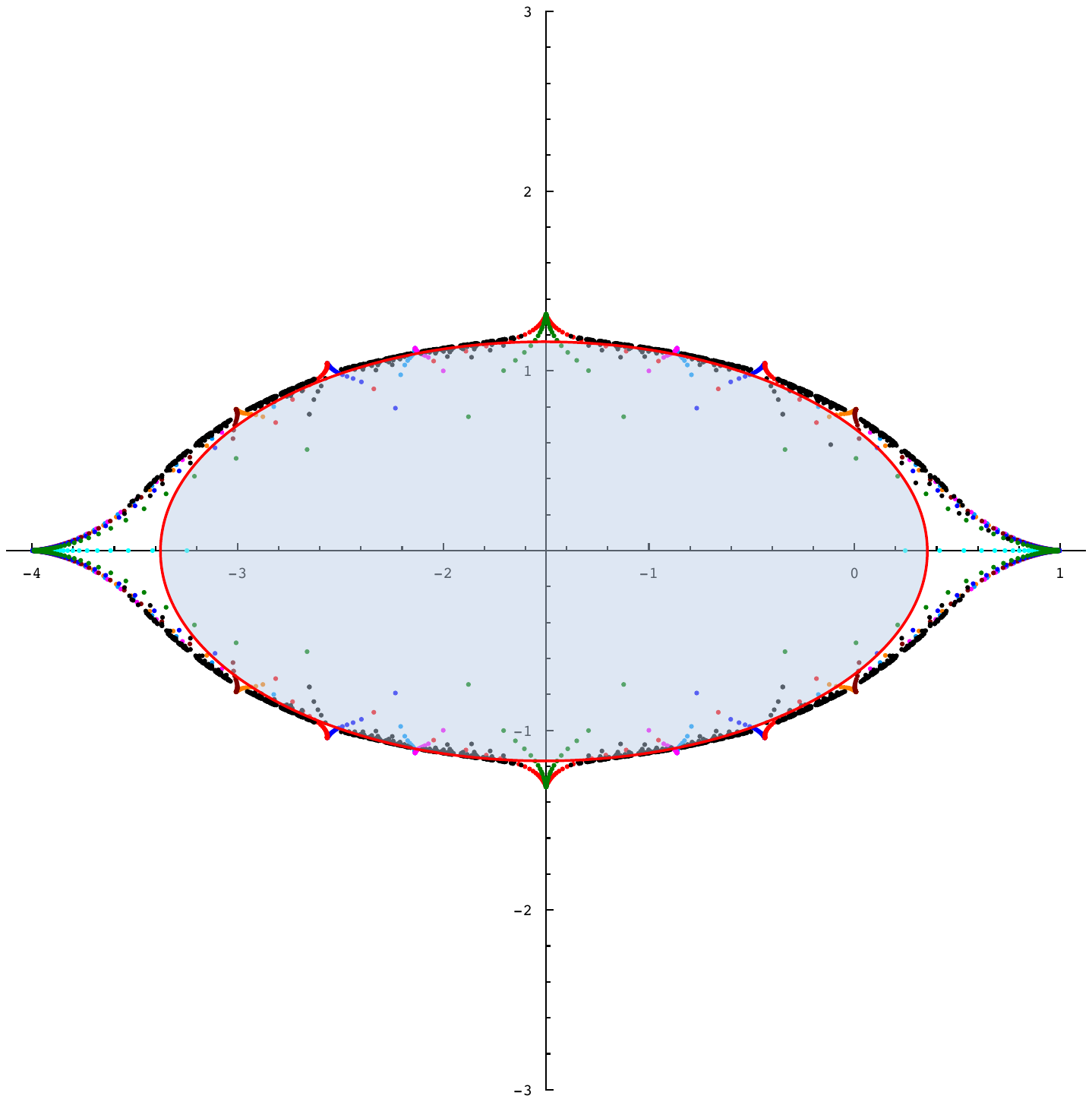} } \\
\end{center}
{\bf Figure 1.}{\em The complement of the bounded region consists of commutator values for discrete groups isomorphic to the free group $\IZ_2*\IZ_3$.  An approximation to the best inscribed ellipse is drawn enclosing commutator values defining the spectrum $\{\delta_n(2,3)\}$}.

\medskip

The issue here  is that although the commutator values $\gamma$ accumulate at every point of the boundary - in particular, and as we will describe,  the pleating rays from the degeneration of a simple closed curve on $\IS^2_{(p,p,q,q)}$ - they seem to accumulate at outward directed cusps,  there might not be infinitely many inside the inscribed ellipse and this seems quite likely if the point of tangency of the ellipse with the boundary of the ``free'' space is a geometrically infinite group (which we conjecture).

\bigskip

\noindent{\bf Proof of Theorem \ref{ads}}  First we observe that calculating the length of the edge of a fundamental polygon for the $(p,q,r)$ triangle groups, $p,q$ both not $2$ and $r=2,3,4\ldots$,   which connects the fixed points of order $p$ and $q$,  gives the infinite family of distinct sequences identified at (\ref{triangle}) and provides the indicated bound.  

Next,  we may suppose that $\langle f_n,g_n \rangle$ is a sequence of discrete groups generated by primitive (rotation angle $\pm \frac{2\pi}{p}$) elliptics of order $p$ and $q$ with hyperbolic distance $ \delta_k$ between their axes,  with $0<\delta_k< M$, $M<\infty$.  The classification of the elementary groups shows that each member of this sequence must be Kleinian.  We may normalise by conjugacy so that the bisector of the common perpendicular between these axes is the point $(0,0,1)\in \IH^3$. The set of all such elliptic elements of order $p$ and $q$ is compact.  By J\o rgensen's algebraic convergence theorem \cite{Jorgensen} there is a subsequence with $f_{n_k}\to f$, $g_{n_k}\to g$,  $f$ and $g$ elliptic of order $p$ and $q$ respectively,  $\langle f_{n_k},g_{n_k}\rangle \to  \langle f,g\rangle$ algebraically,  $\langle f,g\rangle$ is discrete and we have 
\[ \rho(f,g)=\delta =\lim_{k\to\infty}\delta_k \]
 
In fact  J\o rgensen's theorem tells us the map back,  defined by $f\to f_{n_k}$ and $g\to g_{n_k}$,  is an eventual homomorphism.  A nontrivial relator in $\langle f,g\rangle$ yields a polynomial trace condition on $\tr[f,g]-2$ as the traces of $f$ and $g$ are determined.  The roots of this polynomial are isolated and so as $\langle f_{n_k},g_{n_k}\rangle$ inherits this relator if $k$ is large enough, $\tr[f_{n_k},g_{n_k}]-2=\tr[f,g]-2$ for all sufficiently large $k$ and the groups $\langle f_{n_k},g_{n_k}\rangle $ and $\langle f,g\rangle$ are conjugate as the  three complex numbers $\beta(\tilde{f}),\beta(\tilde{g})$ and 
\begin{equation}\label{gammadef}
\gamma(\tilde{f},\tilde{g}) = \tr[\tilde{f},\tilde{g}]-2,
\end{equation}
uniquely determine the group up to conjugacy if $\gamma\neq0$, \cite{GM1}. While $\gamma=0$ implies $\tilde{f},\tilde{g}$ share a fixed point \cite{Beardon} and $\rho(f,g)=0$ and so the group $\langle \tilde{f},\tilde{g} \rangle$ is not Kleinian.

\medskip

 This shows that for fixed $n$, the only limits of groups achieving the values $\delta_n(p,q)$ are groups free on the two generators.  It follows that the spectrum is discrete below $\delta_\infty(p,q)$ and establishes the stated finiteness condition and the assertion on finite multiplicity.  That the space of discrete groups free on the generators (and so isomorphic to $\IZ_p*\IZ_q$) is connected with nonempty interior shows that the spectrum is continuous above $\delta_\infty(p,q)$.  
 
 All that remains to be proved is the estimate $\delta_n(p,q)\leq \delta_n(2,p)$.  We will show that for every group  $\langle f,g\rangle$ which is discrete with $f$ elliptic of order $p$ and $g$ elliptic of order $q$,  there is a group $\langle f,\phi\rangle$ which is also discrete, has $\phi$ elliptic of order two and has closer disjoint axes,  and so is Kleinian. 
 
 Let $\alpha$ be the common perpendicular between the axes of $f$ and $g$ and suppose $\alpha$ meets the axis of $f$ at $x_0\in \IH^3$.  Then $\alpha, g(\alpha)$ and the geodesic from $x_0$ to $g(x_0)$ form a hyperbolic triangle with two edges ($\alpha,g(\alpha)$) of length $\rho(f,g)$ with vertex on the axis of $g$ subtending the angle $2\pi/q$.  Now $f$ and $gfg^{-1}$ are both elliptics of order $p$ with disjoint axes.  There is an involution $\phi$ whose axis passes orthogonally through the bisector of these two axes so that  $\phi f\phi^{-1}=gfg^{-1}$.  By construction
 \[ \rho(f,\phi)=\frac{1}{2} \rho(f,gfg^{-1}) \leq \frac{1}{2} \rho_{\IH^3}(x_0,g(x_0)) \]
From the law of hyperbolic cosines we have
\[\cosh 2\rho(f,\phi) \leq \cosh^2\rho(f,g)-\sinh^2\rho(f,g)\cos\frac{2\pi}{q} \] which implies
\[\cosh^2\rho(f,\phi) \leq  \cosh^2\rho(f,g)-\sinh^2 \rho(f,g)\cos^2\frac{\pi}{q} \leq  \cosh^2\rho(f,g) \]
and equality holds in this last inequality (and throughout) only when $q=2$. This establishes the theorem. \hfill $\Box$

\bigskip

The main aim of this article is to substantially improve the upper-bound at (\ref{ub}) with a view to applications in determining the first several values of the spectrum.   

\section{The geometry of commutators.}

If $f$ and $g$ are each elliptic or loxodromic, then we define the complex distance $\delta+i\theta$ between the axes of $f$ and $g$ as $\delta=\rho(f,g)$ and $\theta$ as the angle between the two hyperplanes formed by each of the axes with the common perpendicular.  We have the following identities,  \cite{GM1}.
   \begin{eqnarray} 
\label{gammageom}
4 \gamma(f,g) & = &  \beta(f) \beta(g)  \sinh^2(\delta+i\theta)
\end{eqnarray}
where $\delta+i\theta$ is the complex distance, and 
\begin{eqnarray}
\label{cosh2r}
\cosh(2\delta)&=&
\left|\frac{4\gamma(f,g)}{\beta(f)\beta(g)}+1\right|
+\left|\frac{4\gamma(f,g)}{\beta(f)\beta(g)}\right|.\end{eqnarray}
   We are often concerned with the case where one of the
isometries,  say $g$,  is of order 2, in which case $\beta(g)=-4$, and we obtain the simpler form
\begin{eqnarray}
\label{order2cosh2r}
\cosh(2\delta)&=&
|1-\gamma(f,g)/\beta(f)|+|\gamma(f,g)/\beta(f)|. 
\end{eqnarray}
Notice that for fixed $\beta(f)\in\IC$ and fixed $\cosh(2\delta)$ at (\ref{order2cosh2r}),  the set of possible values for $\gamma(f,g)$  form an ellipse,  while for  fixed $\cosh(2\theta)$ we get hyperbola. Thus $\delta$ and $\theta$ give very appealing geometric orthogonal coordinates on $\IC\setminus [\beta,0]$.

\section{Bounds for $\delta_\infty(p,q)$: The continuous part of the spectrum.}\label{mainresult}

We first describe how to get better bounds for $\delta_\infty(2,p)$ when $p$ is large.   

\bigskip

The set
\[ {\cal R} = \{z\in\IC :\Big\langle \left[\begin{array}{cc} 1 & 1 \\ 0& 1 \end{array}\right],\left[\begin{array}{cc} 1 & 0 \\ z  & 1 \end{array}\right] \Big\rangle \subset SL(2,\IC) \;\; \mbox{is discrete and free} \}\]
is usually called the Riley slice. It is clearly symmetric under complex conjugation and as $\gamma(f,g)=z^2$ we have $-{\cal R}={\cal R}$.  Lyndon and Ullman were the first to make substantial advances in the study of this set, \cite{LU}, though there was earlier work of Sanov (1947), Chang, Jennings \& Ree (1958) and Brenner (1961).  But it was Riley's computational investigations and the remarkable connections to knot theory which really aroused mathematical interest. The important results of  Keen and Series, \cite{KS}, Komori and Series \cite{KoS} and also of Akiyoshi, Sakuma, Wada \& Yamashita \cite{ASWY} have given very precise descriptions of the Riley slice ${\cal R} \subset \IC$ (illustrated in Figure 3) including the fact that the finite boundary of ${\cal R}$ is a topological circle.

\medskip

The complement of the Riley slice  appears to have a maximal inscribed disk of radius greater than $\frac{3}{2}$,  we have illustrated the disk of radius two.  
\begin{center}
\scalebox{0.6}{\includegraphics*[viewport=20 400 750 750]{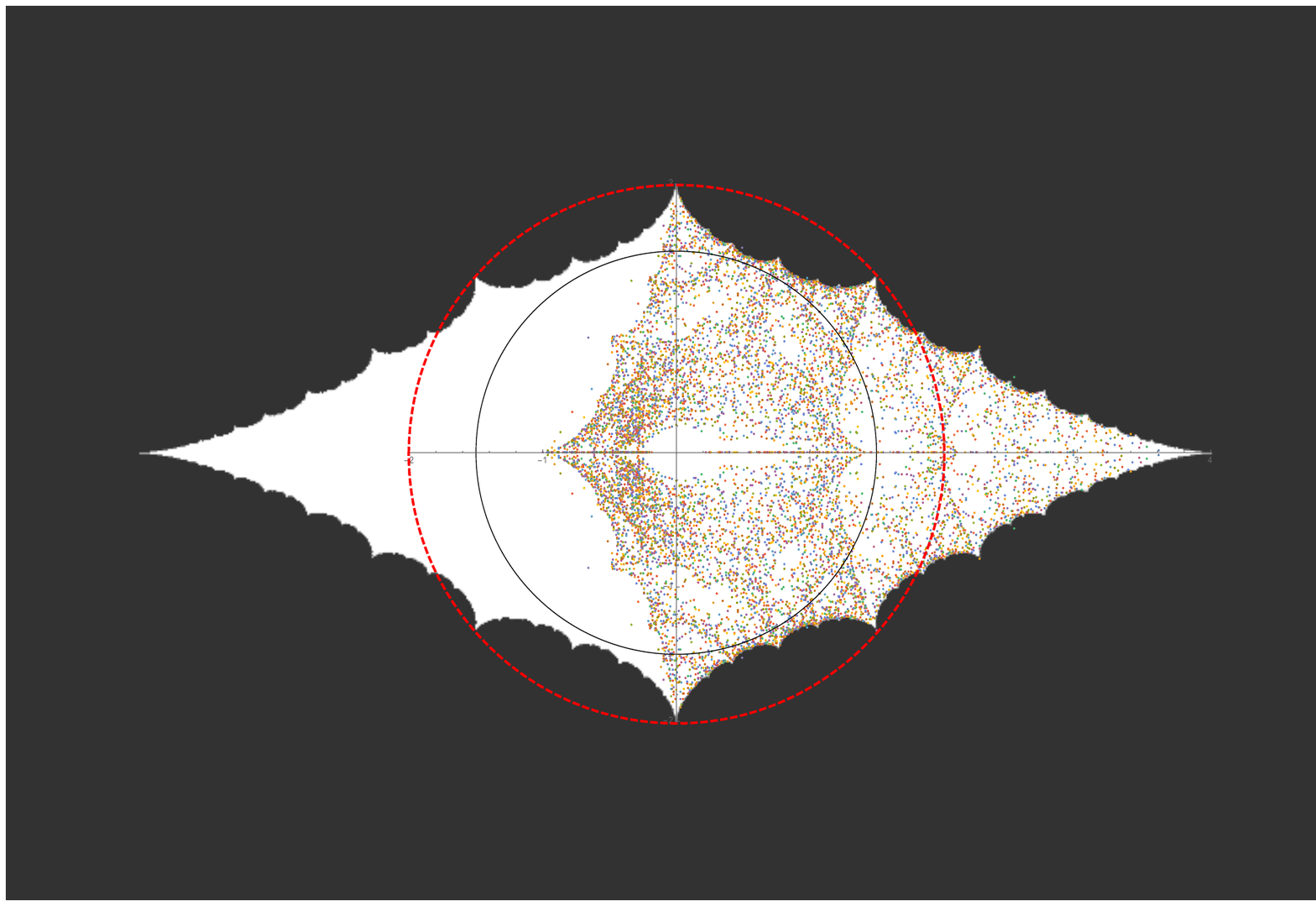} } \\
\end{center}
\noindent{\bf Figure 2.} {\em Roots of the Farey polynomials of the first $257$ Farey sequences and disk of radius $1.5$ and $2$. The Riley slice is the black region.}

\bigskip
To obtain bounds we use the  Keen-Series theory of pleating rays, \cite{KS}.  Briefly, a Farey sequence of order $n$ is the sequence of completely reduced fractions between $0$ and $1$ which when in lowest terms have denominators less than or equal to $n$, arranged in order of increasing size.  Farey sequences enumerate the simple closed geodesics on the four times punctured sphere $\IS^2_4$. Any such surface is obtained by gluing two twice punctured disks together along their common boundary when each has the same length. The conformal structure on $\IS^2_4$ is determined by this length and the twisting angle in $[0,2\pi]$. This outlines why the moduli space of $\IS^2_4$ is an annulus.  The same is true for any $\IS^2_{p,p,q,q}$,  a two sphere with four cone points - two of order $p$ and two of order $q$ each disk enclosing a pair of cone points of the same order.

The hyperbolic length of the geodesic in the homotopy class of the gluing curve is a function of the trace of the word (which as noted represents a simple closed curve representing an element of $\pi_1(\IS^2_4)$) associated with the Farey sequence.  This can be seen through the formula for $\beta(f)$ in terms of the translation length $\tau$ which is the length of the associated geodesic.  As the length of the geodesic on $\IS^2_4$ shrinks to $0$ the surface degenerates to a pair of triply punctured spheres.  This group lies on the boundary of moduli space.  The same holds in the case of four cone points,  the limit being two 2-spheres,  each with one puncture and a pair of cone points - the orbit space of the $(p,p,\infty)$-triangle group.  The relevance to our problem is simply that if $\langle f,g\rangle$ is a geometrically finite Kleinian group free on its two generators of orders $p$ and $q$,  then
\[ \big(\oC \setminus \Lambda(\Gamma)\big)/\Gamma = \IS^2_{p,p,q,q} \]
where $\Lambda(\Gamma)$ is the limit set of $\Gamma$. This is the Schottky representation of this surface group.

\medskip

As Keen and Series recognised, motivated by work of Thurston,  this gives a recipe for finding points on the boundary of the moduli space. Let us work through a very simple example. Consider the fraction $\frac{2}{3}$ given by the Farey sequence  $v=(1,1,-1,1,1-1)$. Choose matrices $X$ and $Y$ as follows
\begin{itemize}
\item In the parabolic case
\[ X_\infty=\left(\begin{array}{cc}1&1\\ 0 &1 \end{array}\right), \quad\quad Y_\infty= \left(\begin{array}{cc}1&0\\ \mu &1 \end{array}\right) \]
\item In the case $f$ has order $p$ and $g$ has order $q$ choose 
 \[ X_p=\left(\begin{array}{cc} e^{i\pi/p}&1\\ 0 &e^{-i\pi/p} \end{array}\right), \quad\quad Y_q= \left(\begin{array}{cc}e^{i\pi/q}&0\\ \mu &e^{-i\pi/q} \end{array}\right) \]
\end{itemize}
Compute the word 
\[ A = X^{v_1}Y^{v_2}X^{v_3}Y^{v_4} \ldots X^{v_{n}}Y^{v_{n}}. \]
In our example of $\frac{2}{3}$ and in the parabolic case,  we compute
\[ A_{\frac{2}{3}}=XYX^{-1}YXY^{-1} = \left(
\begin{array}{cc}
 \mu ^3-\mu ^2+1 & 1-\mu ^2 \\
 (\mu -1)^2 \mu  & -\mu ^2+\mu +1 \\
\end{array}
\right) \]
Then $\tr A_{\frac{2}{3}}+2=p_{2/3}(\mu)=(\mu +1) \left(\mu ^2-3 \mu +4\right)$. The word $A_{\frac{2}{3}}$ will have become parabolic when  $p_{2/3}(\mu)=0$ (so $\tr A_{\frac{2}{3}}=-2$). The three roots of this cubic are $-1$, $\frac{1}{2}(3\pm i \sqrt{7})$.
Keen and Series show that only one complex conjugate pair of roots lies on the boundary of the Riley slice,  the others are either non-discrete groups or discrete groups not free on generators.  In the example we are working through we have the second case here and $\mu=-1$ arises from the Tetrahedral group $Tet(\infty,\infty;3)$,  see \cite{MMec}. Note that $|\frac{1}{2}(3\pm i \sqrt{7})|=2$ and this point lies on the boundary of the Riley slice at the pronounced cusp that meets the circle of radius $2$ illustrated in Figure 2.  

\medskip
We must make the following calculation, recall (\ref{gammadef}).
\begin{equation}
\gamma_{p,q}^\mu=\gamma(A_p,A_q) = \mu_{p,q}  \Big(\mu_{p,q} - 4 \sin \frac{\pi }{p}  \sin \frac{\pi }{q} \Big).
\end{equation}

Actually the collection of all $\gamma_{p,q}^\mu$ where $\mu_{p,q}$ is a root of a Farey polynomials for given $p$ and $q$ must lie in the closure of the complement of the space of freely generated groups.  This gives us a way to illustrate these spaces.  We  include a couple of interesting examples below in Figure 3. In the case of the Riley slice we have $p=q=\infty$ and $\gamma_{\infty,\infty}^\mu=\mu_{\infty,\infty}^2$.  If $f$ and $g$ are both elliptic of the same order, or both parabolic,  then there is an involution $\varphi$ so that $g=\varphi f \varphi$.  It is an elementary calculation to find the trace identity
\[ \gamma(f,g) = \gamma(f,\varphi)(\gamma(f,\varphi)-\beta(f)). \]
Thus,  when $p=q$, we have 
\[ \gamma_{p,p}^\mu = \mu \Big(\mu - 4 \sin^2 \frac{\pi }{p}\Big)=\mu_{p,p}(\mu_{p,p}+\beta) \]
and so 
\begin{equation}
\gamma_{p,2}^\mu =- \mu_{p,p}.
\end{equation}
Thus,  for instance,  the Riley slice can also be identified as the set of commutator values for Kleinian groups $\langle f,\phi\rangle$ freely generated,  where $f$ is parabolic and $\phi$ has order two.

\medskip

In practise it is not easy to sort out the roots of a Farey polynomial to find which indeed lie on the boundary.  This is usually done by tracing out the rays coming from $\infty$ of $p_{r}^{-1}([0,\infty))$.  However that is not a problem we have.  We simply evaluate the formula (\ref{cosh2r}) at each value of $\gamma_{p,q}^{\mu}$.  We can be sure that the ellipse defined by taking the largest value obtained is either meets the boundary (so it is defined by a point on the boundary) or meets the space of freely generated groups as it  encloses the root on the boundary.  Once we identify that
\begin{eqnarray*}
|\gamma_{p,q}^\mu+4\sin^2\frac{\pi}{p}\sin^2\frac{\pi}{q}| & = & |\mu  \Big(\mu - 4 \sin \frac{\pi }{p}  \sin \frac{\pi }{q} \Big)+4\sin^2\frac{\pi}{p}\sin^2\frac{\pi}{q}| \\
& = & |\mu  - 2\sin \frac{\pi }{p}  \sin \frac{\pi }{q}|^2,  
\end{eqnarray*}
together with (\ref{cosh2r}) we have sketched the following theorem.

\begin{theorem} Let $P(\mu)=\tr A_{r/s}^{p,q}(\mu)+2$ associated to the a Farey fraction $r/s$.  Let the roots of $P$ be denoted $\mu_0,\mu_1\ldots\mu_n$ arranged so that  
\[ |\mu_i-2\sin\frac{\pi}{p}\sin\frac{\pi}{q}|^2+|\mu_i||\mu_i-4\sin^2\frac{\pi}{p}\sin^2\frac{\pi}{q}|  \]
is largest for $i=0$. Then
\begin{equation}\label{14}
\cosh 2 \delta_\infty(p,q) \leq \frac{|\mu_0-2\sin\frac{\pi}{p}\sin\frac{\pi}{q}|^2+|\mu_0||\mu_0-4\sin^2\frac{\pi}{p}\sin^2\frac{\pi}{q}| }{4\sin^2 \frac{\pi}{p}\sin^2 \frac{\pi}{p}}
\end{equation}
\end{theorem} 

However we want to get rather more explicit estimates here.  We see that as $p,q\to\infty$,  $X_p\to X_\infty$ and $Y_q\to Y_\infty$. For a fixed Farey rational $r/s$ the polynomials 
$\tr A_{r/s}^{p,q}(\mu) \to \tr A_{r/s}^{\infty,\infty}(\mu)$ locally uniformly in $\IC$.  Hence the roots converge,  and hence the inscribed ellipses will also converge.  Thus to obtain good bounds for $p,q$ large,  we need to first examine the case $p=q=\infty$.  There (\ref{cosh2r}),  appropriately interpreted as $\beta=0$,  asks us to find the largest inscribed disk in the complement of the Riley slice.

\bigskip

Examining these roots of polynomials we settle on the Farey sequence for $\frac{10}{17}$, $v=$  
{\small (1,1,-1,-1 ,1,1, -1, 1,1,-1,-1,1, -1,-1, 1,1, -1, 1,1, -1,-1, 1,1,
-1, 1,1,-1,-1, 1, -1,-1, 1,1,-1)} and then set
\begin{equation}\label{10/17} A = X^{v_1}Y^{v_2}X^{v_3}Y^{v_4} \ldots X^{v_{33}}Y^{v_{34}}
\end{equation}
The trace of $A$ is a polynomial in $\mu$.  In this case it is a degree $17$ polynomial 
$\tr(A)+2=P_{\frac{10}{17}}(\mu)$
\begin{footnotesize}
\begin{eqnarray*} & = & \mu ^{17}-5 \mu ^{16}+18 \mu ^{15}-45 \mu ^{14}+91 \mu ^{13}-151 \mu ^{12}+210 \mu ^{11}-252 \mu ^{10}+255 \mu ^9-225 \mu ^8 \\ && \quad +166 \mu ^7-102 \mu ^6+52 \mu ^5-16 \mu ^4+4 \mu ^3+2 \mu ^2+4 \label{pa}
\end{eqnarray*}
\end{footnotesize}
While not an issue to us as we are just seeking the root of largest modulus, using various polynomial trace identities as discussed in \cite{GM1} we can identify the  root on the boundary of the Riley slice as the other $8$ roots (there is one real root and 7 complex conjugate pairs) are isolated in the space of discrete groups - should they even be discrete. This root is
\[ \mu_0 =0.593302\ldots - i1.42172 \ldots \]
of modulus $|\mu|=1.54055\ldots$.  All these roots are isolated and the next two smaller modulus of roots are $1.4879\ldots$ and $1.3501\ldots$.  We use this information to establish the following Theorem which covers the case $p$, $q$ large.

\begin{theorem}\label{gap} For $p,q\geq 629$,  there is a freely generated Kleinian group $\langle f,g\rangle$,  generated by elements of order $p$ and $q$ on the boundary of the moduli space of such groups and for which  
\[ |\gamma|+|\gamma+4\sin^2\frac{\pi}{p}\sin^2\frac{\pi}{q} | \leq 5, \quad \gamma=\gamma(f,g) \]
Consequently
\begin{equation}
\delta_\infty(p,q) \leq \frac{1}{2}\arccosh\left( \frac{5}{4\sin^2\frac{\pi}{p}\sin^2\frac{\pi}{q}}\right)
\end{equation}
\end{theorem}
\noindent{\bf Remark.} A computational exploration of the Riley slice shows it very unlikely that the constant $5$ can be replaced by any number less than or equal to $4.75$,  we will actually prove a slightly better bound than $5$,  but our approach cannot be sharp as the boundary points we find lie at the end of pleating rays and so on outward directed cusps, \cite{ASWY}.  If we compare with the value $\delta_1(p,q)$ we find for $p\neq q$, 
\begin{eqnarray*} \delta_1(p,q) &=&\arccosh\left( \frac{1}{2\sin\frac{\pi}{p} \sin\frac{\pi}{p}}\right) \approx  \log\;  \frac{pq}{\pi^2} \\
\delta_\infty(p,q) &\approx &\log\;  \frac{\sqrt{5/2} \; pq}{ \pi^2} 
\end{eqnarray*}
which shows a remarkably small gap between the base of the discrete spectrum and the start of the continuous spectrum,
\[ \delta_\infty(p,q)-\delta_1(p,q) \to \frac{1}{2}\log \frac{5}{2}  =0.458145 ,\quad\quad p,q\to\infty. \]
Actually a computational exploration suggests that $0.458145$ is an upper bound for this gap as soon as $p,q\geq 30$,  while of course both terms are growing to $\infty$.

\medskip

\noindent{\bf Proof for Theorem \ref{gap}.} What we first want to do is show that for $p,q$ large,  there is a point on the boundary of moduli space of Kleinian groups freely generated by elements of order $p$ and $q$ within three decimal places of $\mu_0$. We proceed as in the parabolic case with the word $A$ defined at (\ref{10/17}) but with generators
 \[ X_p=\left(\begin{array}{cc} e^{i\pi/p}&1\\ 0 &e^{-i\pi/p} \end{array}\right), \quad\quad Y= \left(\begin{array}{cc}e^{i\pi/q}&0\\ \mu &e^{-i\pi/q} \end{array}\right) \]
The resulting polynomial $P_{\frac{10}{17}}^{p,q}(\mu) $ is too long to write here.  However  is is of degree $17$ in $\mu$ and we can express it as a polynomial in $e^{i\pi/p}-1$ and $e^{i\pi/q}-1$ and multiplying by $e^{i10\pi/p}e^{i10\pi/q} $ will clear the terms with negative exponent.
\[  e^{i10\pi/p}e^{i10\pi/q} P_{\frac{10}{17}}^{p,q}(\mu)  = \sum_{i,j} a_{ij}(\mu) (e^{i\pi/p}-1)^i(e^{i\pi/q}-1)^j \]

 Each of the polynomials $a_{ij}(\mu)$ has integer coefficients and has degree no more than $17$ and $i+j\leq 28$.
 The constant term is the polynomial $P_{\frac{10}{17}} (\mu)$.  We want to show that $P_{\frac{10}{17}}^{p,q}(\mu) $ has roughly the same set of roots if $p,q$ are large.  This is immediate by compactness and the fact the roots depend continuously on the coefficients,  but we would like an explicit estimate.
 
 \begin{lemma} Let $r_0=0.0005$.  Then $|P_{\frac{10}{17}} (\mu_0+re^{i\theta})|\geq 0.0685$.
 \end{lemma}
If $\mu_i$ is any root of  $P_{\frac{10}{17}}$ we compute that  
\[ 28.336  \leq |P_{\frac{10}{17}}'(\mu_i)|\leq139.082\ldots \]
 This shows that $P_{\frac{10}{17}}$ is conformal near each root.  The second derivatives are a bit larger so we proceed by differentiating $|P_{\frac{10}{17}} (\mu_0+re^{i\theta})|^2$ with respect to $\theta$,  computationally finding the two roots and then evaluating the minima (it is about $0.06914005$). \hfill $\Box$
 
 \medskip
 In exactly the same way we have
 \begin{lemma} Let $r=0.025$ and $\mu_i\neq\mu_0$ a roots of $P_{\frac{10}{17}}$.  Then $|P_{\frac{10}{17}} (\mu_i+re^{i\theta})|\geq 0.64$.
 \end{lemma}

Now if $|e^{i\pi/p}-1|,|e^{i\pi/q}-1|< s= 0.005$,  using the same methods we find that 
\begin{equation}
\sum_{i,j} |a_{ij}(\mu_i+re^{i\theta})| s^{i+j} \leq 0.03.
\end{equation}
and 
\begin{equation}
\sum_{i,j} |a_{ij}(\mu_i+re^{i\theta})| s^{i+j} \leq 0.5.
\end{equation}
 
 \begin{corollary} If $|e^{i\pi/p}-1|,|e^{i\pi/q}-1|<\frac{1}{200}$,  then there is a root $\mu_0^{p,q}$ of $P_{\frac{10}{17}}^{p,q}$ with $|\mu_0^{p,q}-\mu_0|<\frac{1}{2000}$. All the other roots of $P_{\frac{10}{17}}^{p,q}$ have modulus smaller than $|\mu_0^{p,q}|$.
 \end{corollary}
 \noindent{\bf Proof.} We have set things up so that on the circles centered at $\mu_i$ and of radius $r=0.025$,  or in the case $i=0$,  $r_0=\frac{1}{2000}$, the following inequality holds
 \[ |P_{\frac{10}{17}}(\mu_i+r_i e^{i\theta})-P_{\frac{10}{17}}^{p,q}(\mu_i+r_i e^{i\theta})| < |P_{\frac{10}{17}}(\mu_i+r_i e^{i\theta})|+|P_{\frac{10}{17}}^{p,q}(\mu_i+r_i e^{i\theta})| \]
 This is so as for $i\neq 0$
\[
|P_{\frac{10}{17}}(\mu_i+r_i e^{i\theta})-P_{\frac{10}{17}}^{p,q}(\mu_i+r_i e^{i\theta})| \leq \sum_{i,j} |a_{ij}(\mu_i +r_i e^{i\theta})|| e^{i\pi/p}-1|^i|e^{i\pi/q}-1|^j  \leq 0.5, 
\]
 and when $i=0$,
\[ |P_{\frac{10}{17}}(\mu_0+r_0 e^{i\theta})-P_{\frac{10}{17}}^{p,q}(\mu_0+r_0 e^{i\theta})| \leq \sum_{i,j} |a_{ij}(\mu_0 +r_0 e^{i\theta})|| e^{i\pi/p}-1|^i|e^{i\pi/q}-1|^j  \leq 0.03.
\]
On the right hand side we have for $i\neq 0$
\[ |P_{\frac{10}{17}}(\mu_i+r_i e^{i\theta})|+|P_{\frac{10}{17}}^{p,q}(\mu_i+r_i e^{i\theta})|>|P_{\frac{10}{17}}^{p,q}(\mu_i+r_i e^{i\theta})| >  0.64 \] 
 and when $i=0$,
\[ |P_{\frac{10}{17}}(\mu_0+r_0 e^{i\theta})|+|P_{\frac{10}{17}}^{p,q}(\mu_0+r_0e^{i\theta})|>|P_{\frac{10}{17}}^{p,q}(\mu_0+r_0 e^{i\theta})| >  0.0685 \] 
The corollary now follows from  the symmetric form of Rouch\'e's theorem.
\begin{lemma}[{\cite[Theorem 3.6 pg 341]{Palka} }]  If $f$ and $g$ are holomorphic and $|f-g|\leq |f|+|g|$ on a simple closed contour,  then $f$ and $g$ have the same number of roots in the region bounded by the contour.
\end{lemma}
We must also observe that the numbers we chose guarantee that $\mu_0^{p,q}$ remains the root of largest modulus. \hfill $\Box$

\medskip

To complete the proof of Theorem \ref{gap} we first observe that 
\[ |e^{i\pi/m}-1|<\frac{1}{200} \Rightarrow m\geq 629. \]
Then we must  use (\ref{order2cosh2r}) and evaluate
\begin{equation}\label{muellipse} |\gamma_{p,q}^\mu+4\sin^2\frac{\pi}{p}\sin^2\frac{\pi}{q}|+|\gamma_{p,q}^\mu| \end{equation}
With $p,q\geq 629$ we find that this number is approximately $2|\gamma_{p,q}^\mu|=2|\mu_0|^2$,  at least to three decimal places as $4\sin^2\frac{\pi}{p}\sin^2\frac{\pi}{q}\leq 10^{-8}$.  Hence
\[ \cosh 2\delta_\infty(p,q) \leq \frac{2|\mu_0|^2}{4\sin^2\frac{\pi}{p}\sin^2\frac{\pi}{q}} = \frac{4.94357}{4\sin^2\frac{\pi}{p}\sin^2\frac{\pi}{q}}  \]
This is enough to establish the theorem.  \hfill $\Box$

\bigskip

\section{$p,q\leq 10$ and $p,\infty$.}

In this section we simply report our computational investigations.  These provide rigorous bounds as described above.  For each pair $(p,q)$ we identify the best Farey rational $r/s$ among the $125$ shortest ($257$ shortest for $p$ small or $\infty$)  and the root of the associated polynomial  $P_{r/s}$ that we were able to find.  Regarding the computation here,  we used Mathematica with root finding algorithms set to an accuracy goal of 40 digits.  We then algebraically found the polynomial in the three variables $\mu$,  $\zeta=e^{i\pi/p}$ and $\xi=e^{i\pi/q}$.  We then evaluated and sorted through all the Farey polynomials and all their roots.  For each polynomial we reported the value obtained at the root lying on the largest ellipse with foci $0$ and $4\sin^2\frac{\pi}{p}\sin^2\frac{\pi}{q}$.   This short circuited having to determine if a root was actually on the boundary of moduli space - note that is is either on the boundary or inside the complement of the free space.  Observationally this seemed to work well and that the roots we obtained appeared from observation (and plotting all the roots)  ``on the boundary''.  It is not the case in general that the roots inside the largest ellipse are not on the boundary.  For example with $p=4$ and $q=\infty$ the polynomial $P_{\frac{21}{25}}(z)$ has the largest ellipse encircling its roots determined by the real root $\mu= 2.59713552133$ which is easily proved to be in the interior of that moduli space, so some other root must lie on the boundary.

\begin{center}
\scalebox{0.3}{\includegraphics*[viewport=20 10 830 930]{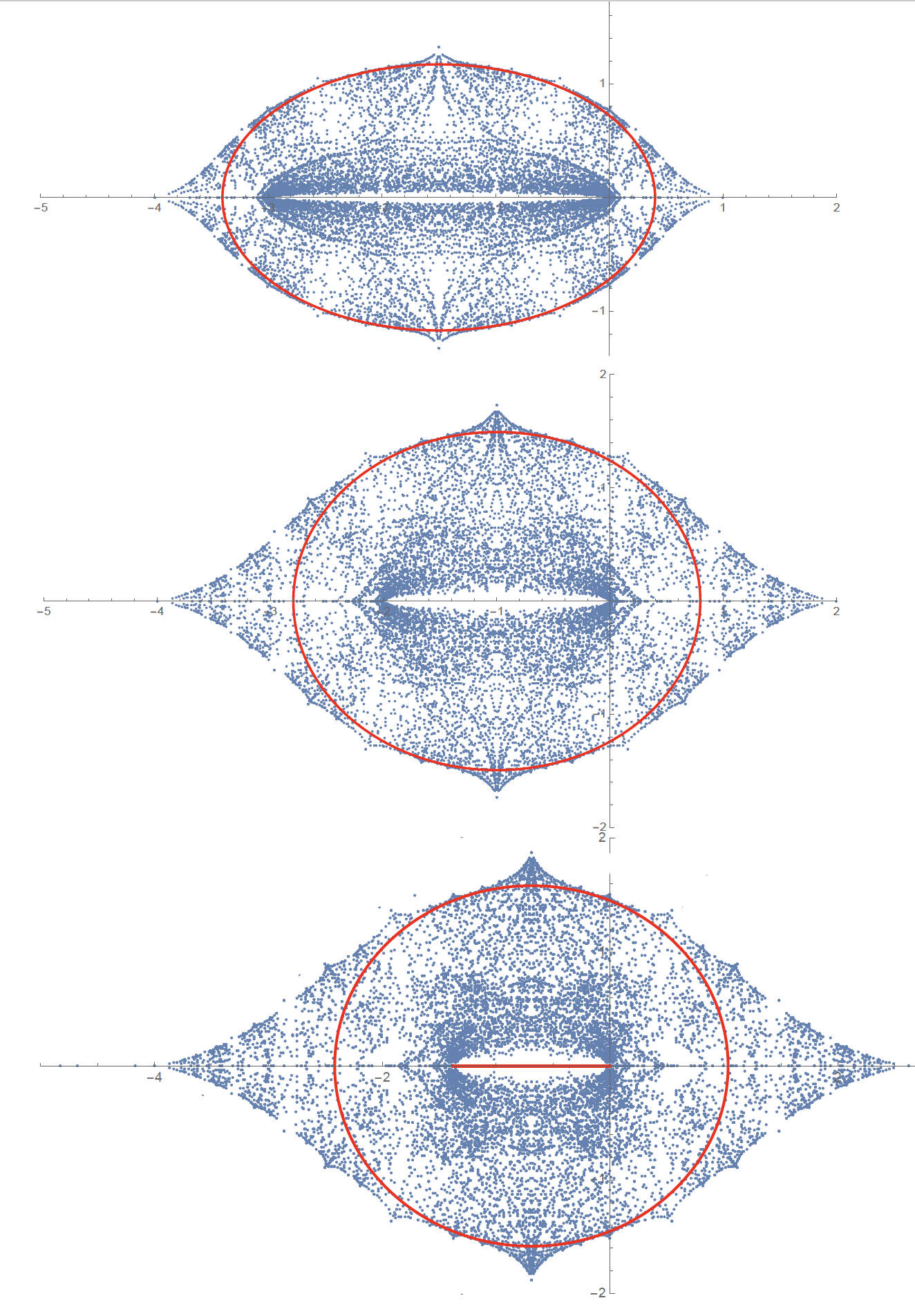} } \\
\end{center}
{\bf Figure 3.}{\em Top: $(2,3)$. Middle: $(2,4)$. Bottom $(2,5)$.  In each case the roots of the Farey polynomials of  degree less than $70$ are illustrated.  The complement of the bounded region they form consists of commutator values for discrete groups isomorphic to the free group $\IZ_2*\IZ_p$,  $p=3,4,5$. The best inscribed ellipse is drawn enclosing commutator values defining the spectrum $\{\delta_n(2,p)\}$}.
\medskip

For applications,  in particular the determination of collaring radii for elliptics,  the most useful bounds are those of $\delta_\infty(2,p)=\frac{1}{2}\delta_\infty(p,p)$ and so these results are presented separately.  The first few cases ($p=3,4,5$) are illustrated above in Figure 3.  The case $p=3$ should also be compared with the Dehn surgery data illustrated in Figure 1.

\bigskip

{\bf Table 1. bounds for $\delta_\infty(2,p)$}\\
\begin{tabular}{|c|c|c|c|c|c|}
\hline
$p$ & Farey Fraction  & upper-bound for $\delta_\infty(2,p)$  & & $\delta_1(2,p)$ \\
 \hline
 $3$ & $ 29/41 $ & $ 0.3583\ldots $  & & $0.19474\ldots$ \\   
\hline
$4$ & $29/41 $ & $ 0.5952 \ldots$ & & $0.4154\ldots $\\
\hline
$5$ & $29/41$ & $0.7831\ldots$ & & $0.4561\ldots$\\
\hline
$6$ & $26/33 $ & $0.9399\ldots$ & & $0.6574\ldots$\\
\hline
$7$ & $ 26/33 $ & $ 1.07948\ldots $  & & $ 0.5452 \ldots $ \\   
\hline
$8$ & $23/29$ & $   1.2013 \ldots $ & & $0.7642 \ldots$\\
\hline
$9$ & $23/29$ & $  1.3105 \ldots $ & & $0.9275 \ldots $\\
\hline
$10$ & $23/29$ & $  1.3979\ldots $  & & $1.0612\ldots$\\
\hline 
$20$ & $ 23/29$ & $ 2.0804 \ldots$  & & $1.8297\ldots $\\
\hline 
$100$ & $ 23/29 $ & $ 3.6821 \ldots $  & & $3.4596\ldots $\\
\hline 
\hline
$p\to\infty$ & $23/29$ & $  {\rm arccosh} \left(\frac{1.2478}{2\sin\frac{\pi}{p}} \right)$ && $ {\rm arccosh}\left(\frac{1}{2\sin\frac{\pi}{p}}\right)$\\
\hline
\end{tabular}

\bigskip
The values identified as $p\to\infty$ are found as follows.  We compute the Farey Polynomial for  $23/29$ with a parabolic generator and an elliptic of order two.  This polynomial of degree $29$ is
\begin{eqnarray*} p_{\frac{23}{29}}(\mu) & = & -\mu ^{29}+5 \mu ^{27}-17 \mu ^{25}+40 \mu ^{23}-74 \mu ^{21}+111 \mu ^{19}-137 \mu ^{17}+144 \mu ^{15} \\&& -126 \mu ^{13}+94 \mu ^{11}-58 \mu ^9+28 \mu ^7-12 \mu ^5+3 \mu ^3-\mu +2 \end{eqnarray*}
and remarkably factors into the product of degree $14$ and $15$ polynomials.  The root of largest modulus is
\[ z_{\frac{23}{29}} = 1.03791 + 0.692732 i, \quad   \big|z_{\frac{23}{29}}\big| = 1.24785. \]
We have $\gamma_{2,\infty}=z_{\frac{23}{29}}^2$ and  inscribed ellipse (recall we now have $\beta=0$ so it is in fact a circle) has 
$ |\gamma|+|\gamma-\beta|=2\big|z_{\frac{23}{29}}\big|^2 = 3.11425\ldots$.
Thus for $p$ large we have 
$\cosh 2\delta_\infty(2,p) \approx \frac{3.11425}{4\sin^2\frac{\pi}{p}}$ 
and hence 
\[   \cosh \delta_\infty(2,p) \approx \sqrt{\frac{3.11425}{4\sin^2\frac{\pi}{p}}+1} \approx \frac{1.2478}{2\sin\frac{\pi}{p}} \]

Notice then that the gap  
\begin{eqnarray*}
\delta_\infty(2,p)-\delta_1(2,p) &=& {\rm arccosh} \left(\frac{1.2478}{2\sin\frac{\pi}{p}} \right)- {\rm arccosh}\left(\frac{1}{2\sin\frac{\pi}{p}}\right)  \\ &&  \approx  \log 1.2478  = 0.2214\ldots
\end{eqnarray*}
independent of $p$.  We will discuss this gap more in the next section.

\bigskip

For more general $p$ and $q$ we have the following theorem.

\begin{theorem} For each $p,q\geq 3$,  the spectrum of maximal collar radii is continuous above the indicated upper-bound.  
\begin{equation}
\cosh\big(2\delta_\infty \big)\leq  \frac{c(p,q)}{4 \sin^2\frac{\pi}{p}\sin^2\frac{\pi}{q}} ,
\end{equation}
where the constants $c(p,q)$ are identified in Table 2. below and the bound is obtained from the Farey polynomial indicated in Table 3. 
\end{theorem}

  {\bf Table 2. Constant $c(p,q)$ }\\
  \begin{footnotesize}
\begin{tabular}{|c|c|c|c|c|c|c|c|c|c|c|c|}
\hline
       & $3$ & $4$ & $5$ & $6$ & $7$ & $8$ & $ 9 $ & $10$ & $20$ & $100$ \\
 \hline
$3$ & $4.984$ & $5.222$ & $5.275$ & $5.279$ & $5.276$ & $ 5.273 $ & $5.268$ & $5.262$ & $5.213$  & $ 5.190$ \\   
\hline
$4$ &   & $5.453$  & $5.512$  & $5.517$ & $5.507$ & $5.496$ & $ 5.487 $ & $5.480$ & $5.453$ & $ 5.443$ \\
\hline
$5$ &   &  & $ 5.484 $ & $5.467$ & $5.451$ & $ 5.440 $ & $5.428$ & $5.433$ & $5.392$ & $5.376$  \\
\hline
$6$ &  & &   & $5.369$ & $5.372$ & $ 5.346$ & $5.327$ & $5.314$ & $5.271$ & $5.258$   \\
 \hline
$7$ &   &  &  &   &  $5.319$ & $5.282$ & $ 5.257 $ & $5.238$ & $5.177$ & $5.157$   \\   
\hline
$8$ &   &   &   &  &  & $ 5.239 $ & $5.209$ & $5.060$ & $5.113$ & $5.089$ \\
\hline
$9$ &   &   &   &  &  &   & $5.173$  & $5.150$ & $ 5.068$ &$ 5.041$\\
\hline
$10$ &   &  &   &  & &  &  &$4.890$ & $5.0352$  &$5.006$  \\
\hline
\\
\hline
$p,q\to \infty$ & $5.16$ & $5.432$ & $5.356,$ & $5.238$  &$5.136$ &  & & & $4.862$ & $$ \\
\hline 
\end{tabular}
\end{footnotesize}

\bigskip

  {\bf Table 3. Farey Polynomial giving $c(p,q)$ }\\
  \begin{footnotesize}
\begin{tabular}{|c|c|c|c|c|c|c|c|c|c|c|c|}
\hline
       & $3$ & $4$ & $5$ & $6$ & $7$ & $8$ & $ 9 $ & $10$ & $20$ & $100$ \\
 \hline
$3$ & $29/41$ & $ {21}/{34} $ & ${19}/{30}$ & ${19}/{30}$ & ${19}/{30}$ & ${19}/{30}$ &${19}/{30}$ & ${19}/{30}$ & $ {16}/{23}  $ & $ {16}/{23} $    \\   
\hline
$4$ &     & $29/41$ & ${18}/{31}$ & ${19}/{31}$ & ${19}/{31}$ & ${19}/{30}$ & $ {19}/{31} $ & ${19}/{30}$ & ${21}/{34}$ & ${21}/{34}$ \\
\hline
$5$ &   &   &$29/41$& $18/31$ & $18/31$ & $18/31 $ & $18/31$ & $18/31$ &  $17/29$ & $19/31$  \\
\hline
$6$ &   &   &  & $26/33$ & $18/31$ & $18/31 $ & $18/31$ & $18/31$ & $18/31$ & $18/31$   \\
 \hline
$7$ &   &   &   &   & $26/33$ & $ 18/31$ & $18/31$ & $18/31$ & $18/31$ &$18/31$  \\   
\hline
$8$ &   &   &   &   &  &  $23/29$ & $18/31$ & $18/31$ & $18/31$ & $18/31$  \\
\hline
$9$ &   &  &   &  &   &   & $23/29$ & $17/29$ & $18/31$ &  $ 18/31 $  \\
\hline
$10$ &  &   & &   &  &   &  & $23/29 $ & $18/31 $  & $18/31 $  \\
\hline
\\
\hline
$ \infty$ & $ 34/49$ & $55/89$ &$46/75$ & $44/75$ & $41/70$ &  & & & $41/70 $ & $23/29$ \\
\hline 
\end{tabular}
\end{footnotesize}

Notice that for each $p\geq 8$, the $(2,3,p)$ triangle group has smallest maximal collar between elliptics of order $p$,  but no other Fuchsian group has maximal collar radius less than $\delta_\infty(2,p)$ and so apart from the smallest, all  Fuchsian maximal collar radii lie in the continuous part of the spectrum.  To see this it is not difficult to see that the next smallest distance between elliptics of order $p$ must among Fuchsian groups occurs in the $(r,p,p)$-triangle groups.  That distance,  for $p\geq 5$ is
\[\cosh  \Delta_2(p,p) \in \frac{\cos\frac{\pi}{r}+\cos^2\frac{\pi}{p}}{\sin^2\frac{\pi}{p}} \geq \cot^2\frac{\pi}{p},
\]
 and this value occurs in the $(2,4,p)$ triangle groups.  Then
 \[ \cosh  \Delta_2(2,p) = \frac{\sqrt{2}}{2\sin\frac{\pi}{p}} \geq \cosh  \delta_\infty(2,p) \]
 for $p\geq 8$ from Table 1.

 \section{An upper bound on the spectral gap.}
 
 Our estimates above show that the gap $\delta_\infty(p,q)-\delta_1(p,q)$ is in general surprisingly small.  Here we give a uniform bound.
 
 \begin{theorem} The maximal spectral gap
 \begin{equation}\label{22} 
 \max_{p,q} \;\; \delta_\infty(p,q) -\delta_1(p,q) \leq \arccosh \;\Delta \approx 1.40598\ldots
 \end{equation}
 \[ \Delta=\frac{1}{3} \Big(\cot ^2 \frac{\pi }{7} +2 \cot  \frac{\pi }{7}  \csc  \frac{\pi }{7} -\sqrt{2 \sin \frac{3 \pi }{14}-1} \left(1+2 \cos  \frac{\pi }{7} \right) \csc ^2 \frac{\pi }{7} \Big) \]
 \end{theorem}
\noindent{\bf Proof.}  We have identified $\delta_1(p,q)$ in Theorem \ref{delta1}.  If $p,q\geq 629$  we may use Theorem\ref{gap} to obtain the bound
\[ \delta_\infty(p,q) - \delta_1(p,q) \leq \frac{1}{2} \arccosh \frac{5}{4\sin^2\frac{\pi}{p}\sin^2\frac{\pi}{q} }  - \arccosh \frac{1}{2\sin\frac{\pi}{p}\sin\frac{\pi}{q} }   <\frac{1}{2} \]
We are then left with a finite problem.  We run through all the bounds for $p,q\leq 629$ using bound of (\ref{6}) obtained from the $(p,q,\infty)$-triangle group  in place of $\delta_\infty$. This value achieves its maximum when $p=3$ and $q=7$ and is the reported value. \hfill $\Box$

\medskip

Actually in the case $p=3$, $q=7$ we can use our better estimate in $\delta_\infty(3,7)$ to get the bound
\[ \frac{1}{2} \arccosh \frac{5.276}{4 \sin ^2 \frac{\pi }{p} \sin ^2 \frac{\pi }{q}} -\arccosh \frac{\cos \frac{\pi }{q} }{2 \sin  \frac{\pi }{p} \sin  \frac{\pi }{q}}  \approx 0.841 \]
This is as large as we were able to find without too much effort.

\begin{theorem}
The bound at (\ref{22}) is sharp for Fuchsian groups and  is achieved for $p=3$, $q=7$.
\end{theorem}
\noindent {\bf Proof.}
For Fuchsian groups generated by elements of order $p$ and $q$ it is straightforward to establish that 
\[ \Delta_1(p,q)=\min_r \arccosh \left[ \frac{\cos\frac{\pi}{r}+\cos\frac{\pi}{p}\cos\frac{\pi}{q}}{\sin\frac{\pi}{p}\sin\frac{\pi}{q}}\right], \quad\quad \frac{1}{r} < 1-\frac{1}{p}-\frac{1}{q},\]
the minimum occuring in the $(p,q,r)$-triangle group.  Also 
\[ \Delta_\infty(p,q)=\arccosh \left[\frac{1+\cos\frac{\pi}{p}\cos\frac{\pi}{q}}{\sin\frac{\pi}{p}\sin\frac{\pi}{q}}\right] \]
the latter occuring in the $(p,q,\infty)$-triangle group. 
If $\frac{1}{p}+\frac{1}{q}<\frac{1}{2}$ we may put $r=2$ to obtain the bound 
\[ \arccosh \left[\frac{1+\cos\frac{\pi}{p}\cos\frac{\pi}{q}}{\sin\frac{\pi}{p}\sin\frac{\pi}{q}}\right] -  \arccosh \left[ \frac{\cos\frac{\pi}{p}\cos\frac{\pi}{q}}{\sin\frac{\pi}{p}\sin\frac{\pi}{q}}\right] .\]
Then we use the difference formula for $\cosh(a+b)$ to see we need to minimise the value of
\begin{eqnarray*}
\Delta & = & \cot \frac{\pi }{p}\cot \frac{\pi }{q}\left(\cot \frac{\pi }{p}\cot  \frac{\pi }{q} +\csc \frac{\pi }{p}\csc \frac{\pi }{q}\right)\\
& & -\sqrt{\cot ^2\frac{\pi }{p}\cot ^2\frac{\pi }{q}-1} \; \left(\csc \frac{\pi }{p}\cot \frac{\pi }{q}+\cot \frac{\pi }{p}\csc \frac{\pi }{q}\right) 
\end{eqnarray*}
Again, it is relatively straightforward to see this occurs at  $p=3$ and $q=7$,  or vice versa.  If $p=2$,  then $q\geq 3$ and we have the bound
\[ \arccosh \left[\frac{1}{\sin\frac{\pi}{q}}\right] - \arccosh \left[ \frac{\cos\frac{\pi}{r}}{\sin\frac{\pi}{q}}\right] , \quad \frac{1}{r}+\frac{1}{q}<\frac{1}{2} \]
Similar bounds are found for $p=3,4$.  This quickly reduces to a finite combinatorial problem and we simply search through all possibilities to see they never achieve the bounds obtained from the $(p,q,2)$ and $(p,q,\infty)$ groups.

\end{document}